\documentclass[letterpaper, 10 pt, conference]{ieeeconf}  

\IEEEoverridecommandlockouts                              

\usepackage{graphicx,cite,algorithmicx,algorithm,algpseudocode,amsthm,amsmath,amsfonts,verbatim,mathtools,enumerate,bm,hyperref,dsfont,tikz,filecontents,pgfplots,subcaption,amssymb}

\usepgfplotslibrary{fillbetween}
\pgfplotsset{compat=1.6}

\allowdisplaybreaks

\newcommand{\diag}{\mathop{\bf diag}}

\newcommand{\argmin}{\mathop{\rm argmin}}

\newcommand{\blkdiag}{\mathop{\rm blkdiag}}
\newcommand{\norm}[1]{\left\lVert#1\right\rVert}
\newcommand{\mnorm}[1]{{\left\vert\kern-0.25ex\left\vert\kern-0.25ex\left\vert #1 
    \right\vert\kern-0.25ex\right\vert\kern-0.25ex\right\vert}}

\newcommand{\cvx}{\texttt{C}}
\newcommand{\ncvx}{\texttt{N}}

\newcommand{\ie}{{\it i.e.}}

\title{\LARGE \bf Filtering-Linearization: A First-Order Method for Nonconvex Trajectory Optimization with Filter-Based Warm-Starting
}

\author{Minsen Yuan, Ryan J. Caverly, and Yue Yu
\thanks{
M. Yuan, R. J. Caverly, and Y. Yu are with the Department of Aerospace Engineering and Mechanics at the University of Minnesota Twin Cities, Minneapolis, MN 55455. (email:\{yuan0450, rcaverly, yuey\}@umn.edu) Y. Yu would like to thank Iman Askari for helpful early discussions.
}   
}

\begin{document}

\maketitle
\thispagestyle{empty}
\pagestyle{empty}

\begin{abstract}
Nonconvex trajectory optimization is at the core of designing trajectories for complex autonomous systems. A challenge for nonconvex trajectory optimization methods, such as sequential convex programming, is to find an effective warm-starting point to approximate the nonconvex optimization with a sequence of convex ones. We introduce a first-order method with filter-based warm-starting for nonconvex trajectory optimization. The idea is to first generate sampled trajectories using constraint-aware particle filtering, which solves the problem as an estimation problem. We then identify different locally optimal trajectories through agglomerative hierarchical clustering. Finally, we choose the best locally optimal trajectory to warm-start the prox-linear method, a first-order method with guaranteed convergence. We demonstrate the proposed method on a multi-agent trajectory optimization problem with linear dynamics and nonconvex collision avoidance. Compared with sequential quadratic programming and interior-point method, the proposed method reduces the objective function value by up to approximately 96\% within the same amount of time for a two-agent problem, and 98\% for a six-agent problem. 

\end{abstract}

\pgfplotsset{every axis plot}
\pgfplotsset{grid style=dotted}

\section{Introduction}

One fundamental problem in optimal control is to optimize the trajectories of dynamical systems subject to various constraints. The objective of this optimization is to minimize a cost function that typically evaluates reference tracking performance and control effort~\cite{malyuta2021advances,malyuta2021convex}. The constraints include physical constraints, such as bounds on kinematics, and operational constraints, such as collision avoidance. Trajectory optimization plays a critical role in the control of a wide range of autonomous systems, including aircraft, spacecraft, and ground vehicles~\cite{malyuta2021advances,malyuta2021convex,wang2024survey}.
Trajectory optimization becomes computational demanding when nonconvex constraints are present. For example,  collision avoidance often requires the trajectories remain within a nonconvex set, leading to nonconvex constraints. These nonconvex constraints make the problem computationally expensive to solve, especially in large-scale problems.

A common approach to solving trajectory optimization with nonconvex constraints is to construct a sequence of subproblems based on a starting trajectory, each of which approximates the original optimization problem. For example, the subproblems in the Interior-Point Method (IPM)~\cite{el1996formulation,byrd1999interior,lukvsan2004interior} approximate the inequality constraints in the original problem using log barrier functions. Sequential Quadratic Programming (SQP)~\cite{boggs1995sequential,gill2011sequential} solves a sequence of quadratic programming subproblems, each of which approximates the nonlinear problem with a quadratic program. Sequential Convex Programming (SCP)~\cite{morgan2014model,vitus2008tunnel} uses a sequence similar to SQP but replaces the quadratic program with a convex problem. A common challenge in these existing methods is the lack of a starting trajectory sufficiently close to the optimum. As a result, these methods often suffer from slow convergence, or convergence to suboptimal or even infeasible trajectories.


A recent method, constraint-aware particle filtering~\cite{askari2021nonlinear,askari2022sampling,askari2023model}, can avoid the need for warm-starting. This approach reformulates a nonconvex trajectory optimization problem as a nonlinear state estimation problem, using particle filtering to approximate the distributions of the optimal trajectories. However, this method typically requires a large number of particles to achieve highly accurate trajectories, especially in large-scale problems, leading to heavy computational demands. Additionally, this method does not clearly address the problems where multiple locally-optimal trajectories exist, as it provides a distribution of all locally-optimal trajectories without indicating which one is globally optimal. 


We propose a first-order method with filter-based warm-starting for trajectory optimization with linear dynamics and nonconvex inequality constraints. This method is a novel combination of the prox-linear method---an SCP method with guaranteed convergence~\cite{drusvyatskiy2018error,drusvyatskiy2019efficiency}---and constraint-aware particle filtering. The key contribution is that we use particle filtering to warm-start the prox-linear method with a small number of particles to find an optimal trajectory for the original nonconvex problem with reasonable computational effort. First, we generate a set of sampled trajectories using constraint-aware particle filtering that approximate the distribution of optimal trajectories. Second, we identify multiple locally-optimal trajectories by applying hierarchical clustering \cite{nielsen2016hierarchical} to the sampled trajectories. Third, we select the best center of the cluster based on the lowest cost value and minimal constraint violations. Finally, we use this best center as the warm-starting trajectory for initializing the prox-linear method.

We demonstrate the proposed method on a multi-agent trajectory optimization problem, where each agent must avoid collision with static obstacles as well as other agents. To compare our proposed method with the benchmark methods, namely IPM and SQP, in terms of the convergence of the objective function value and constraints violation over time, we conduct 100 Monte Carlo simulations to evaluate their performance. The proposed approach demonstrates both faster convergence and lower trajectory costs compared to benchmark methods. It reduces the objective function value by up to approximately 96\% compared to SQP and IPM with random initialization within the same amount of time for a two-agent problem, and by 98\% for a six-agent problem.


\section{Prox-Linear Method for Nonconvex Trajectory Optimization}

We introduce the implementation of the prox-linear method~\cite{drusvyatskiy2018error,drusvyatskiy2019efficiency} for trajectory optimization with linear dynamics and nonconvex inequality constraints.

\subsection{Trajectory Optimization with Nonconvex Constraints}
We consider the following trajectory optimization problem
\begin{equation}\label{opt: traj}
    \begin{array}{ll}
   \underset{x_{1:N}, u_{1:N}}{\mbox{minimize}}  & \sum_{k=1}^N \norm{Cx_k-\hat{y}_k}_Q^2+ \norm{u_k}_R^2 \\
    \mbox{subject to} & x_1=\hat{x}_1,\enskip x_{k+1}=Ax_k+Bu_k, \\
    & g(x_k, u_k)\leq 0_{n_g}, \enskip 1\leq k \leq N,
    \end{array}
\end{equation}
where \(Q\in\mathbb{S}_{\succ 0}^{n_y}\), \(R\in\mathbb{S}_{\succ 0}^{n_u}\) are positive definite cost matrices, \(C\in\mathbb{R}^{n_x}\to\mathbb{R}^{n_y}\) is a matrix that maps state to regulated output. Vector \(\hat{x}_1\in\mathbb{R}^{n_x}\) is the initial state. Vector \(\hat{y}_k\in\mathbb{R}^{n_x}\) is a reference state for all \(k=1, \ldots, N\); \(A\in\mathbb{R}^{n_x\times n_x}\) and \(B\in\mathbb{R}^{n_x\times n_u}\) are dynamics parameters; \(g:\mathbb{R}^{n_x}\times \mathbb{R}^{n_u}\to\mathbb{R}^{n_g}\) is a continuously differentiable function describing the state and input constraints.

\subsection{Prox-Linear Method}

To pinpoint the nonconvex inequality constraints in optimization~\eqref{opt: traj}---which are the sources of difficulty in solving it---we introduce the following partition of function \(g\)
\begin{equation}
    g(x_k, u_k)\coloneqq \begin{bmatrix}
        g_\cvx(x_k, u_k)\\
        g_\ncvx(x_k, u_k)
    \end{bmatrix},
\end{equation}
where \(g_\cvx:\mathbb{R}^{n_x}\times \mathbb{R}^{n_u}\to \mathbb{R}^{n_\cvx}\), \(g_\ncvx:\mathbb{R}^{n_x}\times \mathbb{R}^{n_u}\to \mathbb{R}^{n_\ncvx}\), with \(n_\cvx+n_\ncvx=n_g\) and the set \(\mathbb{G}\coloneqq \{x_k\in\mathbb{R}^{n_x}, u_k\in\mathbb{R}^{n_u}|g_\cvx(x_k, u_k)\leq 0_{n_\cvx}\}\) is a convex set. With this partition, we can re-write optimization~\eqref{opt: traj} equivalently as
\begin{equation}\label{opt: traj cnc}
    \begin{array}{ll}
   \underset{x_{1:N}, u_{1:N}}{\mbox{minimize}}  & \sum_{k=1}^N \norm{Cx_k-\hat{y}_k}_Q^2+ \norm{u_k}_R^2 \\
    \mbox{subject to} & x_1=\hat{x}_1,\enskip  x_{k+1}=Ax_k+Bu_k, \\
&        g_\cvx(x_k, u_k)\leq 0_{n_\cvx}, \enskip 
        g_\ncvx(x_k, u_k) \leq 0_{n_\ncvx},\\
        & 1\leq k\leq N.
    \end{array}
\end{equation}

The prox-linear method solves optimization~\eqref{opt: traj cnc} via iteratively linearizing the nonlinear dynamics constraints and nonconvex inequality constraints \cite{drusvyatskiy2019efficiency,drusvyatskiy2018error}. To introduce this algorithm, we first introduce the following \emph{penalized optimization}, which approximates the nonconvex constraints in optimization~\eqref{opt: traj} via penalty functions
\begin{equation}\label{opt: penalty 2}
    \begin{array}{ll}
   \underset{\substack{x_{1:N}, u_{1:N},\\ r_{1:N}}}{\mbox{minimize}}  & \sum_{k=1}^N \big(\norm{Cx_k-\hat{y}_k}_Q^2+ \norm{u_k}_R^2 +\gamma \mathbf{1}_{n_\ncvx}^\top s_k\big)  \\
    \mbox{subject to} & x_1=\hat{x}_1,\enskip x_{k+1}=Ax_k+Bu_k, \\
    & g_\cvx(x_k, u_k)\leq 0_{n_\cvx},\\
    & g_\ncvx (x_k, u_k)\leq s_k,\enskip s_k\geq 0_{n_\ncvx},\enskip 1\leq k\leq N.
    \end{array}
\end{equation}
Let \(\{x_{1:N}^j, u_{1:N}^j\}\) denote the solution at the \(k\)-th iteration (\(k\in\mathbb{N}\)) of prox-linear method, and
\begin{equation}
\begin{aligned}
    & H_k^j \coloneqq \partial_x g_\ncvx(x_k^j, u_k^j), \enskip L_k^j \coloneqq \partial_u g_\ncvx(x_k^j,u_k^j),\\
    &d_k^j \coloneqq g_\ncvx(x_k^j,u_k^j) -\partial_x g_\ncvx(x_k^j,u_k^j) x_k^j-\partial_u g_\ncvx(x_k^j,u_k^j) u_k^j,
\end{aligned}
\end{equation}
where \(\gamma\in\mathbb{R}_{>0}\) is a large penalty weight. The \((k+1)\)-th iteration of the prox-linear method solves an approximated version of \eqref{opt: penalty 2} as\begin{equation}\label{opt: linear}
    \!\!\!\begin{array}{ll}
    \underset{\substack{x_{1:N}, u_{1:N},\\ r_{1:N}}}{\mbox{minimize}}  & \sum_{k=1}^N \big(\norm{Cx_k-\hat{y}_k}_Q^2+ \norm{u_k}_R^2+\gamma \mathbf{1}_{n_\ncvx}^\top s_k\big) \\
   & +\frac{1}{2\rho}\sum_{k=1}^N \big(\norm{x_k-x_k^j}^2_2+ \norm{u_k-u_k^j}^2_2\big)\\
    \mbox{subject to} &x_1=\hat{x}_1,\enskip g_\cvx(x_k, u_k)\leq 0_{n_\cvx},\\
    &x_{k+1}=A x_k+B u_k,\\
    &H_k^j x_k+L_k^ju_k+d_k^j\leq s_k,\enskip s_k\geq 0_{n_\ncvx},\\
    & 1\leq k\leq N,
    \end{array}
\end{equation}
where \(\rho>0\) is a step size. We summarize the complete prox-linear method in Algorithm~\ref{alg: prox-linear}. The idea is to linearize the constraints in \eqref{opt: traj} around the trajectories obtained in the previous iteration, then penalize the deviation from the linearization point in the objective function.

\begin{algorithm}
\caption{Prox-Linear Method}
\label{alg: prox-linear}
\begin{algorithmic}[1]
\Require Initial trajectory \(\{x_{1:N}^0, u_{1:N}^0\}\). Constant \(\rho\in(0, 1)\), \(\epsilon, \delta\in\mathbb{R}_{> 0}\). Let \(j=0\) and \(s_k=\epsilon\mathbf{1}_{n_\ncvx}\). 
\While{\(\displaystyle\sum_{k=1}^N\big(\norm{x_k^j-\overline{x}_k^j}^2_2+\norm{u_k^j-\overline{u}_k^j}^2_2\big)> \epsilon\) or \(\displaystyle\sum_{k=1}^N\norm{s_k}_2^2 > \epsilon\)}
\State Solve optimization~\eqref{opt: linear} with \(\rho=1/\gamma\), obtain optimal trajectory \(\{\overline{x}_{1:N}^j, \overline{u}_{1:N-1}^j\}\)
\State \(\{x_{1:N}^{j+1}, u_{1:N}^{j+1}\}\gets \{\overline{x}_{1:N}^j, \overline{u}_{1:N}^j\}\)
\State \(j\gets j+1\)
\EndWhile
\end{algorithmic}
\end{algorithm}

The convergence of Algorithm~\eqref{alg: prox-linear} hinges on the initial trajectory \(\{x_{1:N}^0, u_{1:N}^0\}\). Its key idea is to repeatedly linearize nonconvex constraints. The closer it starts from the optimum, the more accurate the linearization, and faster the convergence. 
\section{Warm-Starting via Constraint-Aware Particle Filtering and Clustering}

The challenge in using Algorithm~\eqref{alg: prox-linear} is to find an initial trajectory \(\{x_{1:N}^0, u_{1:N}^0\}\) sufficiently close to the optimum. As our main contribution, we propose a method to compute warm-starting trajectories by combining constraint-aware particle filtering~\cite{askari2021nonlinear,askari2022sampling,askari2023model} and hierarchical clustering~\cite{nielsen2016hierarchical}.

\subsection{Constraint-Aware Particle Filtering}

The constraint-aware particle filtering method solves optimization~\eqref{opt: traj} by solving an moving-horizon estimation problem via sampling \cite{askari2021nonlinear,askari2022sampling, askari2023model}. To this end, we first approximate optimization~\eqref{opt: traj} as follows:

\begin{equation}\label{opt: estimation}
    \begin{array}{ll}
    \underset{\substack{x_{1:N}, u_{1:N}\\
    w_{0:N-1}, v_{1:N}, q_{1:N}}}{\mbox{minimize}} &  \sum_{k=1}^N \big(\norm{v_k}_Q^2+\norm{w_{t-1}}_R^2+  \norm{q_k}^2 \big)\\
    \mbox{subject to} & \begin{bmatrix}
        x_1\\
        u_1
    \end{bmatrix} =\begin{bmatrix}
        \hat{x}_1\\
        w_0
    \end{bmatrix},\\
    & \begin{bmatrix}
        x_{k+1}\\
        u_{k+1}
    \end{bmatrix}=\begin{bmatrix}
        Ax_k+Bu_k\\
        0_{n_u}
    \end{bmatrix}+\begin{bmatrix}
        0_{n_x}\\
        w_k
    \end{bmatrix},\\
    & \begin{bmatrix}
        \hat{y}_k\\
        -\nu \mathbf{1}_{n_g}
    \end{bmatrix}=\begin{bmatrix}
        Cx_k\\
        [g(x_k, u_k)]_+
    \end{bmatrix}+\begin{bmatrix}
        v_k\\
        q_k
    \end{bmatrix},\\
    &1\leq k\leq N,
    \end{array}
\end{equation}
where \([z]_+\coloneqq \max\{z, 0_{n_g}\}\) for all \(z\in\mathbb{R}^{n_g}\), \(\nu\in\mathbb{R}_{>0}\) is a weighting parameter. 

By choosing an appropriate value for \(\nu\), we can obtain an optimal solution for optimization~\eqref{opt: traj} by solving its penalized variant~\eqref{opt: estimation}. In particular, we can rewrite the objective function in~\eqref{opt: estimation} equivalently as
\begin{equation*}
    \sum_{k=1}^N \left(\norm{Cx_k-\hat{y}_k}_Q^2+\norm{u_k}_R^2+\norm{[g(x_k, u_k)]_++\nu \mathbf{1}_{n_g} }_2^2\right).
\end{equation*}
In other words, optimization~\eqref{opt: estimation} replaces the inequality constraints in optimization~\eqref{opt: traj} with a linear and quadratic penalty for their violations. The exact penalty theorem states that if \(\{x_{1:N}, u_{1:N}, w_{0:N-1}, v_{1:N}, q_{1:N}\}\) is optimal for optimization~\eqref{opt: estimation} with a sufficiently large \(\nu\) and \(\{x_{1:N}, u_{1:N}\}\) satisfies the constraints in optimization~\eqref{opt: traj}, then \(\{x_{1:N}, u_{1:N}\}\) is also optimal for optimization~\eqref{opt: traj} \cite[Thm.~17.4]{nocedal1999numerical}.

\begin{algorithm}
\caption{Unscented Transform}
\label{alg: unscented}
\begin{algorithmic}[1]
\Require Mean \(x\in\mathbb{R}^n\), variance \( A_1\in\mathbb{S}_{\succ 0}^{n}, A_2\in\mathbb{S}_{\succeq 0}^l\), function \(\varphi:\mathbb{R}^n\to\mathbb{R}^l\), \(\theta \approx 0.1\).
\State \(\lambda\gets (\theta^2-1)n\)
\State \(\overline{X}\gets x\mathbf{1}_{2n+1}^\top\)
\State \(X\gets \overline{X}+\sqrt{n+\lambda}\begin{bmatrix}
0_n & \sqrt{A_1} & -\sqrt{A_1}
\end{bmatrix}\)
\State \(a\gets\frac{1}{n+\lambda}\begin{bmatrix}
  \lambda & \frac{1}{2}\mathbf{1}_{2n}^\top  
\end{bmatrix}^\top\)
\State \(b \gets \frac{1}{n+\lambda}\begin{bmatrix}
    \lambda +(n+\lambda)(3-\theta^2) & \frac{1}{2}1^\top_{2n}
\end{bmatrix}^\top\)
\State \(Y\gets\varphi(X)\) (\(\varphi\) applied column-wise)
\State \(y\gets Ya\)
\State\(\overline{Y}\gets y\mathbf{1}_{2n+1}^\top\)
\State \(B_1 \gets (Y-\overline{Y})\diag(b)(Y-\overline{Y})^\top+A_2\)
\State \(B_2 \gets (X-\overline{X})\diag(b) (Y-\overline{Y})^\top\)
\Ensure \(y, B_1, B_2\)
\end{algorithmic}
\end{algorithm}

The constraint-aware particle filtering method considers the following \emph{virtual} linear stochastic system with a nonlinear output function:
\begin{subequations}
    \begin{align}
        \begin{bmatrix}
        x_{k+1}\\
        u_{k+1}
    \end{bmatrix} & = \begin{bmatrix}
        A & B\\
        0_{n_u\times n_x} & 0_{n_u\times n_u}
    \end{bmatrix}\begin{bmatrix}
        x_k\\
        u_k
    \end{bmatrix}+\begin{bmatrix}
        0_{n_x}\\
        w_k
    \end{bmatrix},\\
    \eta_k & =\begin{bmatrix}
        Cx_k\\
        [g(x_k, u_k)]_+
    \end{bmatrix}+\begin{bmatrix}
        v_k\\
        q_k
    \end{bmatrix}, \\
    w_k &\sim\mathcal{N}(0_{n_u}, R^{-1}),\\ v_k&\sim\mathcal{N}(0_y, Q^{-1}),\enskip q_k\sim\mathcal{N}(0_{n_g}, I_{n_g}).
    \end{align}
\end{subequations}
Let
\begin{equation}
\begin{aligned}
    &\xi_k  \coloneqq \begin{bmatrix}
        x_k\\
        u_k
    \end{bmatrix},\, \hat{\eta}_k  \coloneqq \begin{bmatrix}
        \hat{y}_k\\
        -\nu\mathbf{1}_{n_g}
    \end{bmatrix}, \,\overline{A}\coloneqq \begin{bmatrix}
        A & B\\
        0_{n_u\times n_x} & 0_{n_u\times n_u}
    \end{bmatrix},\\ & \psi(\xi_k)\coloneqq\psi\left(\begin{bmatrix}
        x_k\\
        u_k
    \end{bmatrix}\right)\coloneqq \begin{bmatrix}
        Cx_k\\
        [g(x_k, u_k)]_+
    \end{bmatrix},\\
    E & \coloneqq \blkdiag(0_{n_x\times n_x}, R^{-1}), \enskip F \coloneqq \blkdiag(Q^{-1}, I_{n_g}).
\end{aligned}
\end{equation}

The idea of constraint-aware particle filtering is to approximate the maximum likelihood distribution of \(\xi_{1:N}\) using a finite number of samples~\cite{askari2023model}. We summarize the constraint-aware particle filtering method in Algorithm~\ref{alg: filter}, which uses the unscented transform, summarized in Algorithm~\ref{alg: unscented}, as a sub-procedure. Its core idea is to solve optimization~\eqref{opt: estimation} via sequential sampling. In practice, we often approximate the nonsmooth ReLU function with the smooth \emph{softplus function}, \ie, \([g(x_k, u_k)]_+\approx \ln(1+\exp(g(x_k, u_k)))\).

With a small number of particles, Algorithm~\ref{alg: filter} can solve optimization~\eqref{opt: estimation} to suboptimality. The resulting trajectory usually satisfies the constraints in \eqref{opt: traj} while having close to optimal cost. Further improving the trajectory quality requires increasing the number of particles in Algorithm~\ref{alg: filter}, which tends to lead to a diminishing return. However, a suboptimal trajectory generated by Algorithm~\ref{alg: filter} with a small number of particles serves as a good warm-starting point for Algorithm~\ref{alg: prox-linear}: the latter converges rapidly when initialized close to an optimal trajectory \cite{drusvyatskiy2018error}.

\begin{algorithm}
\caption{Constraint-Aware Particle Filtering}
\label{alg: filter}
\begin{algorithmic}[1]
\Require Initial state \(\xi_0^i\in\mathbb{R}^{n_x+n_u}\), variance \(\Sigma_0^i\in\mathbb{S}_{\succeq 0}^{n_x+n_u}\) for all \(i=1, 2, \ldots, m\). Desired output \(\hat{\eta}_k\in\mathbb{R}^{n_x+n_g}\) for all \(k=1, 1, \ldots, N\). Sampling parameter \(\alpha\in(0, 1)\) and \(\kappa\in(1, m)\).
\State Let \(\omega^i_k=\frac{1}{m}\) for all  \(1\leq i\leq m\) and \(1\leq k\leq N.\)
\For{\(k=1, 2, \ldots, N-1\)}
\For{\(i=1, 2, \ldots, m\)}
\State \((\zeta, U, V)\gets\text{Alg.~\ref{alg: unscented}}(\overline{A}\xi_k^i, \overline{A}\Sigma_k^i \overline{A}^\top+E, F, \psi)\) 
\State \(K\gets V (U)^{-1}\)
\State \(\Sigma^i\gets \overline{A}\Sigma_k^i \overline{A}^\top+E-KU K^\top\) 
\State \(z\sim \mathcal{N}(0_{n_x+n_u}, \alpha I_{n_x+n_u})\) 
\State \(\xi_{k+1}^i\gets \overline{A}\xi_k^i+K(\hat{\eta}_{k+1}-\zeta)+\sqrt{\Sigma^i}z\) 
\State \(\tilde{\omega}_{k+1}^i\gets\frac{\omega_k^i}{\sqrt{\det U}} \exp\left(-\frac{1}{2}\norm{\hat{\eta}_{k+1}-\zeta}_{(U)^{-1}}^2\right)\)
\EndFor
\State \(\omega_{k+1}^i \gets \tilde{\omega}_{k+1}^i/\left(\sum_{j=1}^m\tilde{\omega}_{k+1}^j\right)\) for all \(1\leq i\leq m\).
\If{\(\kappa\sum_{j=1}^m(\omega_{k+1}^j)^2\geq 1\) } 
\State Let \(\{\overline{\xi}_s^i\}_{s=0}^{t+1}\gets\{\xi_s^i\}_{s=0}^{t+1}\) and \(\overline{\Sigma}^i\gets\Sigma^i\) for all \(1\leq i\leq m\).
\For{\(i=1, 2, \ldots, m\)}
\State Sample a random integer \(j\) such that \(\mathds{P}(j=l)=\omega_{k+1}^l\) for all \(1\leq l\leq m\).
\State \(\{\xi_s^i\}_{s=0}^{k+1}\gets\{\overline{\xi}_s^j\}_{s=0}^{k+1},\Sigma^i\gets \overline{\Sigma}^j,\omega_{k+1}^i\gets \frac{1}{m}.\)
\EndFor 
\EndIf
\EndFor

\Ensure \(\{\{\xi_k^i, \omega_k^i\}_{k=1}^N\}_{i=1}^m\).
\end{algorithmic}
\end{algorithm}

\subsection{Trajectory Clustering}
The outcome of Algorithm~\ref{alg: filter} provides multiple sampled trajectories, namely \(\{\xi_k^1\}_{k=1}^N\), \(\{\xi_k^2\}_{k=1}^N\), \(\ldots ,\) \(\{\xi_k^m\}_{k=1}^N\). If the trajectory optimization in \eqref{opt: traj} admits multiple optimal solutions, then these trajectories will form multiple \emph{clusters}, each of which contains several trajectories close to one of the local optimum. To identify these clusters, we propose to use \emph{hierarchical clustering}. To this end, we first need to define a \emph{distance matrix} \(D\in\mathbb{R}^{m\times m}\), such that
\begin{equation}
    D_{ij}=\begin{cases}
        0, & \text{if } i=j,\\
        \sum_{k=1}^N  p(\xi_k^i, \xi_k^j), & \text{otherwise.} 
    \end{cases}
    \label{eqn:distance}
\end{equation}
where \(p(\xi_k^i, \xi_k^j)\) measure the distance between vector \(\xi_k^i\) and \(\xi_k^j\). A naive choice for function \(p(\xi_k^i, \xi_k^j)\) is simply the Euclidean norm between \(\xi_k^i\) and \(\xi_k^j\). Here we propose the following quadratic distance based on the cost parameters in optimization~\eqref{opt: traj}
\begin{equation}
    p(\xi_k^i, \xi_k^j)\coloneqq (\xi_k^i-\xi_k^j)^\top\blkdiag (C^\top QC, R)(\xi_k^i-\xi_k^j).
\end{equation}
After calculating the distances between each pair of trajectories, we use agglomerative hierarchical clustering~\cite[Sec.~8.2]{nielsen2016hierarchical} with group average linkage and choose a cutting height to distinguish all the clusters.



\subsection{Cluster Selection}

After obtaining all the clusters, we proceed to select the best cluster through three steps. First, for each cluster, we calculate the averaged trajectory by considering the weighted contributions of the trajectories within the cluster. Let \( n_y \) be the total number of clusters, and \( \mathcal{Y}^j \) represent the set of trajectories in cluster \( j \), where \( 1 \leq j \leq n_y \). For each time step \( k \), the cluster center state \( \bar{\xi}_k^j \) for cluster \( j \) is computed as a weighted sum of the states from all trajectories in this cluster at time step \( k \):
\begin{equation}
\bar{\xi}_k^j = \sum_{i \in \mathcal{Y}^j} \frac{\omega_k^i}{\sum_{i \in \mathcal{Y}^j} \omega_k^i} \xi_k^i,
\end{equation}
where \( \omega_k^i \) is a scalar from Algorithm~\ref{alg: filter}, representing the weight of the \( i \)-th trajectory in the set \( \mathcal{Y}^j \) at time step \( k \). \( \xi_k^i \) denotes the state of the \( i \)-th trajectory at time step \( k \) in the set \( \mathcal{Y}^j \). This process is repeated for all \( N \) time steps to obtain the complete averaged trajectory.

Second, we evaluate each averaged trajectory by summing its objective function value from Algorithm~\ref{opt: traj} and the weighted constraint violations, including violations of the state dynamics constraints and the constraints function \( g \). We let
\begin{equation}
\begin{aligned}
&\textstyle \phi^j\coloneqq \sum_{k=1}^N \left( \norm{C\bar{x}_k^j - \hat{y}_k}_Q^2 + \norm{\bar{u}_k^j}_R^2\right)\\
&\textstyle +\alpha  \norm{ \bar{x}_1^j - \hat{x}_1 }_1  + \alpha\sum_{k=1}^{N-1} \norm{ A\bar{x}_k^j + B\bar{u}_k^j - \bar{x}_{k+1}^j }_1 \\ 
&\textstyle +\alpha\sum_{k=1}^N\norm{\max\{ g(\bar{x}_k^j, \bar{u}_k^j), 0_{n_g} \}}_1.
\end{aligned}
\end{equation}
For all \( j = 1, 2, \ldots, n_y \), where \( \bar{x}_k^j \) and \( \bar{u}_k^j \) are derived from the averaged trajectory \( \bar{\xi}_k^j \), with \(\bar{\xi}_k^j \coloneqq \begin{bmatrix}
    (\bar{x}_k^j)^\top & (\bar{u}_k^j)^\top
\end{bmatrix}^\top\), and where \( \alpha \) is a weighting parameter.

Finally, we warm-start Algorithm~\ref{alg: prox-linear} by letting
\begin{equation}
    \{x_k^0, u_k^0\}_{k=1}^N= \{\bar{x}_k^{j^\star}, \bar{u}_k^{j^\star}\}_{k=1}^N,
    \label{eq: best_center}
\end{equation}
where \(j^\star\in\argmin_{1\leq j\leq n_y} \enskip \phi^j\).


\section{Numerical Experiments on Multiagent Trajectory Optimization}

We demonstrate the application of the proposed method in multi-agent trajectory optimization problems and compare its convergence performance with benchmark methods.

\subsection{Multiagent Trajectory Optimization with Collision Avoidance Constraints}

We consider the special case of optimization~\eqref{opt: traj}, where the state \( x_k \) contains the position and velocity of \(n \in \mathbb{N}\) agents. The trajectory reference \( \hat{y}_k \) is defined as the linear interpolation between the initial and final positions. The input \( u_k \) represents the acceleration. The state dynamics are derived by discretizing the continuous-time model using a Zero Order Hold approach. The matrix \( C \) maps the state \( x_k \) to the position output. The matrices \( A \), \( B \), and \( C \) are defined as
\begin{subequations}
\begin{align}
    &A = \exp(A_c \Delta t),\enskip B = \int_0^{\Delta t} \exp(A_c s) \, ds \cdot B_c,\\
    &C = I_n \otimes \begin{bsmallmatrix} 1 & 0 & 0 & 0 \\ 0 & 1 & 0 & 0 \end{bsmallmatrix},
\end{align} 
\end{subequations}
where
\(A_c = I_n \otimes \begin{bsmallmatrix} 0 & 0 & 1 & 0 \\ 0 & 0 & 0 & 1 \\ 0 & 0 & 0 & 0 \\ 0 & 0 & 0 & 0 \end{bsmallmatrix}\), \(B_c = I_n \otimes \begin{bsmallmatrix} 0 & 0 \\ 0 & 0 \\ 1 & 0 \\ 0 & 1 \end{bsmallmatrix}
\), and \( I_n \) is the identity matrix of size \( n \times n \), \( \otimes \) denotes the Kronecker product, and \( \Delta t \) is the discrete time interval.

We consider the case where the convex constraints \( g_\cvx(x_k, u_k) \leq 0_{n_\cvx} \) denote the following elementwise upper and lower bounds for the velocity and acceleration:
\begin{subequations}
\begin{align}
    &\begin{bmatrix}
        I_n\otimes \begin{bsmallmatrix} 0 & 0 & 1 & 0 \\ 0 & 0 & 0 & 1 \end{bsmallmatrix}\\
        -I_n\otimes \begin{bsmallmatrix} 0 & 0 & 1 & 0 \\ 0 & 0 & 0 & 1 \end{bsmallmatrix}
    \end{bmatrix}x_k - 
    \begin{bmatrix}
        v_{\max}\\
        -v_{\min}
    \end{bmatrix}\leq 0_{4n}, \\
    &\begin{bmatrix}
        u_k\\
        -u_k
    \end{bmatrix} - 
    \begin{bmatrix}
        a_{\max}\\
        -a_{\min} 
    \end{bmatrix} \leq 0_{4n},
\end{align}
\end{subequations}
where \( v_{\min}, v_{\max} \in \mathbb{R}^{2n} \) are the lower and upper bounds for the velocity; and \( a_{\min}, a_{\max} \in \mathbb{R}^{2n} \) are the lower and upper bounds for the acceleration, respectively.

We consider the case where the nonconvex constraints \( g_\ncvx(x_k, u_k) \leq 0_{n_\ncvx}\) include 
\begin{subequations}
\begin{align}
&1- \norm{Z(\theta^l)^\top (M_i x_k - c_e^l)}_{\diag(p^l)}^2\leq 0,\label{eqn: ellipse}\\
&\gamma^2-\norm{M_i x_k - M_j x_k}_2^2 \leq 0,
\label{agents_self}
\end{align}
\end{subequations}
for all \(1\leq i, j\leq n\in\mathbb{N}_{>0}\) such that \(i\neq j\) and \(l=1, 2, \ldots, n_o\) where \(n_o\in\mathbb{N}_{>0}\) is the total number of elliptical obstacles, \(
Z(\theta^l) = \begin{bmatrix} \cos(\theta^l) & -\sin(\theta^l) \\ \sin(\theta^l) & \cos(\theta^l) \end{bmatrix}
\), and \( M_i \) represents the transformation from state \( x_k \) to the position of the \( i \)-th agent. The constraint in \eqref{eqn: ellipse} specifies the elliptical obstacle avoidance constraints. Here, \(\theta^l\in[0, 2\pi]\), \(c_e^l\in\mathbb{R}^2\), and \(p^l\in\mathbb{R}^2\) are the parameters defining the orientation, center position, and semi-axes length of the \(l\)-th elliptical obstacle, respectively. The constraint in \eqref{agents_self} ensures that agents maintain a minimum distance to avoid collisions, where \( M_i \) and \( M_j \) represent the transformations from state \( x_k \) to the positions of the \( i \)-th and \( j \)-th agents, respectively, and \( \gamma \) represents the minimum distance threshold between any two agents.


\subsection{Numerical Experiments}

We consider the cases of \( n = 3 \) and \( n = 6 \), with \( n_o = 3 \), where the parameters of obstacles and the upper and lower bounds for velocity and acceleration are given by
\begin{align*}
c_e^1 &= \begin{bmatrix} 3 & 3 \end{bmatrix}^\top, \enskip p^1 = \frac{1}{4}\begin{bmatrix} 1 & 1 \end{bmatrix}^\top, \enskip \theta^1 = 0, \\
c_e^2 &= \begin{bmatrix} 9 & 5 \end{bmatrix}^\top, \enskip p^2 = \frac{1}{2.25}\begin{bmatrix} 1 & 1 \end{bmatrix}^\top, \enskip \theta^2 = 0, \\
c_e^3 &= \begin{bmatrix} 5 & 9 \end{bmatrix}^\top, \enskip p^3 = \begin{bmatrix} \frac{1}{4} & \frac{1}{2.25} \end{bmatrix}^\top, \enskip \theta^3 = \frac{\pi}{3}, \\
v_{\max} &= - v_{\min} = 2\cdot\mathbf{1}_{2n}, \enskip a_{\max} = - a_{\min} = \mathbf{1}_{2n}.
\end{align*}

\begin{figure}[!t]
\begin{subfigure}{0.49\columnwidth}
\centering
\includegraphics[width=\linewidth]{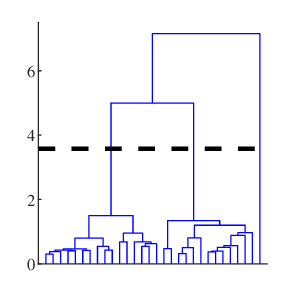}
\caption{Two agents.}
\label{fig:dendrogram_2}
\end{subfigure}
\begin{subfigure}{0.49\columnwidth}
\centering
\includegraphics[width=\linewidth]{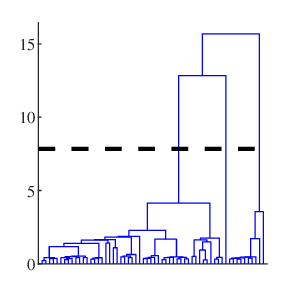}
\caption{Six agents.}
\label{fig:dendrogram_6}
\end{subfigure}
\caption{The dendrogram and cutting height for two and six agents. We illustrate the cutting height using a dashed line to cut the dendrogram and obtain flat partitions. In our case, we set the height to \( 50\% \) of the maximum group average linkage distance.
}
\label{fig:dendrogram}
\end{figure}

We generate sampled trajectories using Algorithm~\ref{alg: filter}, then cluster these trajectories using hierarchical clustering and illustrate the corresponding dendrogram in Fig.~\ref{fig:dendrogram}. Specifically, Fig.~\ref{fig:dendrogram_2} is for two agents with $m = 30$, $\kappa = 12$, and $\Sigma_0^i = I_{n_x+n_u}$, while Fig.~\ref{fig:dendrogram_6} is for six agents with $m = 60$, $\kappa = 24$, and $\Sigma_0^i = \operatorname{diag}([10^{-1}\mathbf{1}_{n_x}; \mathbf{1}_{n_u}])$. The parameters $N = 30$, $\alpha = 5 \times 10^{-3}$, and $\xi_0^i = 0_{n_x+n_u}$ are the same for both cases as those used in Algorithm~\ref{alg: filter}. The leaves at the bottom of these diagrams represent the sampled trajectories. The \(y\)-axis label represents the group average linkage of the branches, with the distance between two sample trajectories calculated using~(\ref{eqn:distance}). We can clearly identify different clusters by choosing the appropriate cutting height in the dendrogram.

\begin{figure}[!htp]
\begin{subfigure}{0.49\columnwidth}
\centering
\pgfplotsset{every tick label/.append style={font=\scriptsize}}

\begin{tikzpicture}
\begin{axis}[ymin=-2, ymax=14,
    xmin=-2, xmax=14,  width=1.2\columnwidth,
    height=1.2\columnwidth,
        xtick = {0, 5, 10},
        ytick = {0, 5, 10}]

\draw[rotate around={-90+60:(axis cs:5,9)},fill=gray!10] (axis cs: 5,9) ellipse (1.5 and 2);

\draw[fill=gray!10] (axis cs:3,3) circle (2);
\draw[fill=gray!10] (axis cs:9,5) circle (1.5);

\addplot[dashed,color=red] table[x=x1, y=y1] {2agent_optimal.dat};

\addplot[dashed,color = blue] table[x=x2, y=y2] {2agent_optimal.dat};

\addplot[solid, thick, color=red] table[x=x1_true, y=y1_true] {2agent_optimal.dat};

\addplot[color=red, mark=*, mark options={fill=red, draw=red, scale=1.2}] coordinates {(0, 0)};
\addplot[color=red, mark=pentagon*, mark options={fill=red, draw=red, scale=1.5}] coordinates {(5.0043, 5.6046)};
\addplot[color=red, mark=triangle*, mark options={fill=red, draw=red, scale=1.5}] coordinates {(12, 12)};

\addplot[solid, thick, color=blue] table[x=x2_true, y=y2_true] {2agent_optimal.dat};

\addplot[color=blue, mark=*, mark options={fill=blue, draw=blue, scale=1.2}] coordinates {(5, 12)};
\addplot[color=blue, mark=pentagon*, mark options={fill=blue, draw=blue, scale=1.5}] coordinates {(6.9642,  6.0029)};
\addplot[color=blue, mark=triangle*, mark options={fill=blue, draw=blue, scale=1.5}] coordinates {(6, 0)};

\end{axis}
\end{tikzpicture}
\caption{Two agents.}
\label{fig: 2agent traj}
\end{subfigure}
\begin{subfigure}{0.49\columnwidth}
\centering
\pgfplotsset{every tick label/.append style={font=\scriptsize}}

\begin{tikzpicture}

\begin{axis}[ymin=-2, ymax=14,
    xmin=-2, xmax=14,  width=1.2\columnwidth,
    height=1.2\columnwidth,
        xtick = {0, 5, 10},
        ytick = {0, 5, 10}]

\draw[rotate around={-90+60:(axis cs:5,9)},fill=gray!10] (axis cs: 5,9) ellipse (1.5 and 2);

\draw[fill=gray!10] (axis cs:3,3) circle (2);
\draw[fill=gray!10] (axis cs:9,5) circle (1.5);

\definecolor{color1}{rgb}{1,0,0} 
\definecolor{color2}{rgb}{0,0,1} 
\definecolor{color3}{rgb}{0,1,0} 
\definecolor{color4}{rgb}{1,0,1} 
\definecolor{color5}{rgb}{0,1,1} 
\definecolor{color6}{rgb}{0.6350,0.0780,0.1840}

\addplot[dashed,color=red] table[x=x1, y=y1] {6agent_optimal.dat};

\addplot[solid, thick, color=red] table[x=x1_true, y=y1_true] {6agent_optimal.dat};

\addplot[color=red, mark=*, mark options={fill=red, draw=red, scale=1.2}] coordinates {(6.5, 13)};
\addplot[color=red, mark=pentagon*, mark options={fill=red, draw=red, scale=1.5}] coordinates {(6.94036, 9.4364)};
\addplot[color=red, mark=triangle*, mark options={fill=red, draw=red, scale=1.5}] coordinates {(6.5, 0.003)};


\addplot[dashed,color = blue] table[x=x2, y=y2] {6agent_optimal.dat};

\addplot[solid, thick, color=blue] table[x=x2_true, y=y2_true] {6agent_optimal.dat};

\addplot[color=blue, mark=*, mark options={fill=blue, draw=blue, scale=1.2}] coordinates {(3.25, 0)};
\addplot[color=blue, mark=pentagon*, mark options={fill=blue, draw=blue, scale=1.5}] coordinates {(6.3339,  3.6479)};
\addplot[color=blue, mark=triangle*, mark options={fill=blue, draw=blue, scale=1.5}] coordinates {(9.7499, 13)};


\addplot[dashed,color = green] table[x=x3, y=y3] {6agent_optimal.dat};

\addplot[solid, thick, color=green] table[x=x3_true, y=y3_true] {6agent_optimal.dat};

\addplot[color=green, mark=*, mark options={fill=green, draw=green, scale=1.2}] coordinates {(9.75, 0)};
\addplot[color=green, mark=pentagon*, mark options={fill=green, draw=green, scale=1.5}] coordinates {(4.8688,  7.1018)};
\addplot[color=green, mark=triangle*, mark options={fill=green, draw=green, scale=1.5}] coordinates {(3.2446, 13.0018)};

\addplot[dashed,color = magenta] table[x=x4, y=y4] {6agent_optimal.dat};

\addplot[solid, thick, color=magenta] table[x=x4_true, y=y4_true] {6agent_optimal.dat};

\addplot[color=magenta, mark=*, mark options={fill=magenta, draw=magenta, scale=1.2}] coordinates {(0, 3.25)};
\addplot[color=magenta, mark=pentagon*, mark options={fill=magenta, draw=magenta, scale=1.5}] coordinates {(7.2163,  5.4741)};
\addplot[color=magenta, mark=triangle*, mark options={fill=magenta, draw=magenta, scale=1.5}] coordinates {(13, 9.75)};

\addplot[dashed,color = cyan] table[x=x5, y=y5] {6agent_optimal.dat};

\addplot[solid, thick, color=cyan] table[x=x5_true, y=y5_true] {6agent_optimal.dat};

\addplot[color=cyan, mark=*, mark options={fill=cyan, draw=cyan, scale=1.2}] coordinates {(0, 9.75)};
\addplot[color=cyan, mark=pentagon*, mark options={fill=cyan, draw=cyan, scale=1.5}] coordinates {(4.4591,  5.1442)};
\addplot[color=cyan, mark=triangle*, mark options={fill=cyan, draw=cyan, scale=1.5}] coordinates {(13, 3.2402)};

\addplot[dashed,color = color6] table[x=x6, y=y6] {6agent_optimal.dat};

\addplot[solid, thick, color=color6] table[x=x6_true, y=y6_true] {6agent_optimal.dat};

\addplot[color=color6, mark=*, mark options={fill=color6, draw=color6, scale=1.2}] coordinates {(13, 6.5)};
\addplot[color=color6, mark=pentagon*, mark options={fill=color6, draw=color6, scale=1.5}] coordinates {(6.8421,  7.4388)};
\addplot[color=color6, mark=triangle*, mark options={fill=color6, draw=color6, scale=1.5}] coordinates {(0, 6.5)};

\end{axis}
\end{tikzpicture}
\caption{Six agents.}
\label{fig: 6agent traj}
\end{subfigure}

\caption{Warm-starting and optimal position trajectories of two and six agents. The dashed lines and solid lines represent the warm-starting and optimal position trajectories, respectively. All trajectories start from the initial point (circle), pass through the midpoint (pentagon), and end at the final point (triangle).}
\label{optplot:trajectory}
\end{figure}
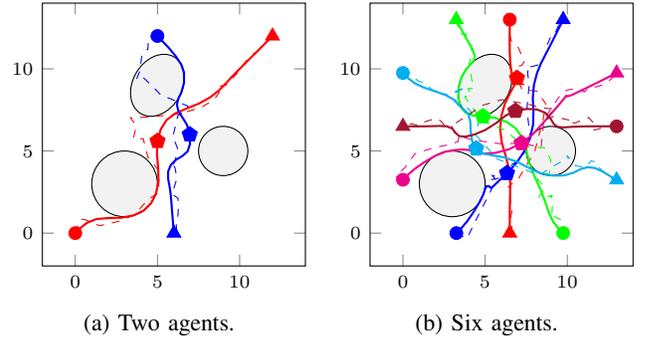


We illustrate the warm-starting and optimal trajectories in Fig.~\ref{optplot:trajectory}. The dashed lines represent the warm-starting trajectories computed according to~\eqref{eq: best_center}, while the solid lines represent the optimal trajectories computed by Algorithm~\ref{alg: prox-linear} using the warm-start in~\eqref{eq: best_center}. Figures~\ref{fig: 2agent traj} and~\ref{fig: 6agent traj} show that the warm-starting trajectories are already close to the optimal trajectories, though they may violate the collision avoidance constraints.

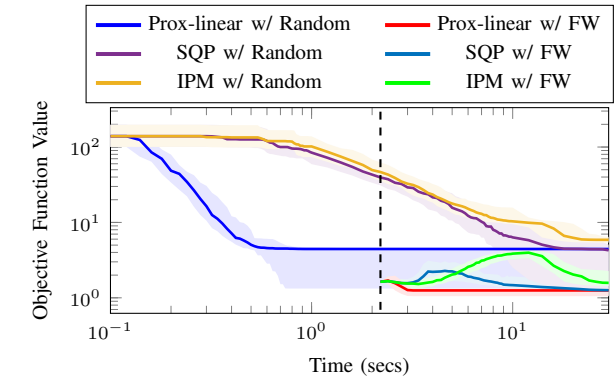
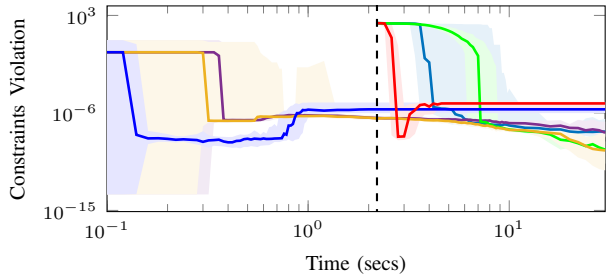
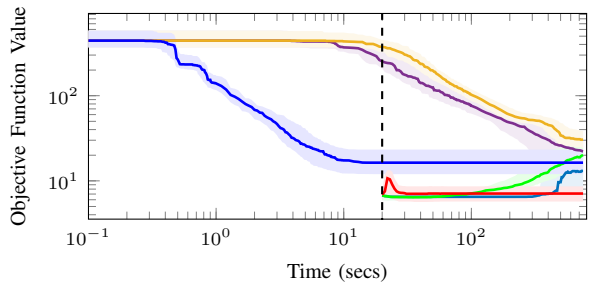
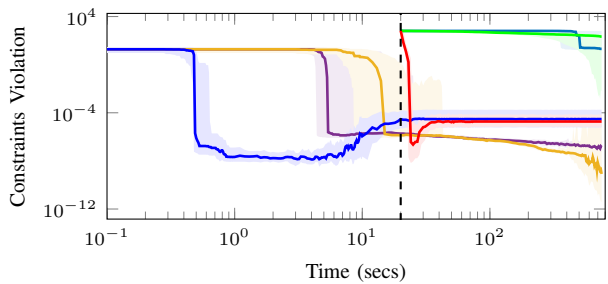
\begin{figure}[!htp]

\begin{subfigure}{\linewidth}
\pgfplotsset{every tick label/.append style={font=\scriptsize}}
\begin{tikzpicture}
\begin{loglogaxis}[
    xlabel = \footnotesize Time (secs), ylabel=\footnotesize Objective Function Value,
    width=0.95\linewidth,
    height = 0.5\linewidth,
    xmin = 1e-1, xmax = 30, xshift = 0.3cm,
    legend style={
        at={(-0.05, 1.5)}, 
        anchor=north west,
        legend columns=2,
        /tikz/column 2/.style={
            column sep=10pt},
        font=\footnotesize
    },
    legend image post style={line width=1.5pt}, 
]


\definecolor{proxrandom_color}{rgb}{0,0,1}
\addplot [mark=none, color=proxrandom_color,line width=1pt] table[x=averageTimes,y=objMedianVals] {2agent_proxrandom.dat};
\addlegendentry{Prox-linear w/ Random} 

\definecolor{proxhc_color}{rgb}{1,0,0}
\addplot [mark=none, color=proxhc_color,line width=1pt] table[x=averageTimes,y=objMedianVals] {2agent_proxhc.dat};
\addlegendentry{Prox-linear w/ FW} 

\definecolor{sqprandom_color}{rgb}{0.4940,0.1840,0.5560}
\addplot [mark=none, color=sqprandom_color,line width=1pt] table[x=averageTimes,y=objMedianVals] {2agent_sqprandom.dat};
\addlegendentry{SQP w/ Random} 

\definecolor{sqphc_color}{rgb}{0,0.4470,0.7410}
\addplot [mark=none, color=sqphc_color,line width=1pt] table[x=averageTimes,y=objMedianVals] {2agent_sqphc.dat};
\addlegendentry{SQP w/ FW} 

\definecolor{ipmrandom_color}{rgb}{0.9290,0.6940,0.1250}
\addplot [mark=none, color=ipmrandom_color,line width=1pt] table[x=averageTimes,y=objMedianVals] {2agent_ipmrandom.dat};
\addlegendentry{IPM w/ Random} 

\definecolor{ipmhc_color}{rgb}{0, 1, 0}
\addplot [mark=none, color=ipmhc_color,line width=1pt] table[x=averageTimes,y=objMedianVals] {2agent_ipmhc.dat};
\addlegendentry{IPM w/ FW} 



\addplot [name path=upper, draw=none] table[x=averageTimes,y=objUpperMedianVals] {2agent_proxrandom.dat};
\addplot [name path=lower, draw=none] table[x=averageTimes,y=objLowerMedianVals] {2agent_proxrandom.dat};
\addplot [fill=proxrandom_color!10] fill between[of=upper and lower];

\addplot [name path=upper, draw=none] table[x=averageTimes,y=objUpperMedianVals] {2agent_sqphc.dat};
\addplot [name path=lower, draw=none] table[x=averageTimes,y=objLowerMedianVals] {2agent_sqphc.dat};
\addplot [fill=sqphc_color!10] fill between[of=upper and lower];

\addplot [name path=upper, draw=none] table[x=averageTimes,y=objUpperMedianVals] {2agent_sqprandom.dat};
\addplot [name path=lower, draw=none] table[x=averageTimes,y=objLowerMedianVals] {2agent_sqprandom.dat};
\addplot [fill=sqprandom_color!10] fill between[of=upper and lower];


\addplot [name path=upper, draw=none] table[x=averageTimes,y=objUpperMedianVals] {2agent_ipmhc.dat};
\addplot [name path=lower, draw=none] table[x=averageTimes,y=objLowerMedianVals] {2agent_ipmhc.dat};
\addplot [fill=ipmhc_color!10] fill between[of=upper and lower];


\addplot [name path=upper, draw=none] table[x=averageTimes,y=objUpperMedianVals] {2agent_ipmrandom.dat};
\addplot [name path=lower, draw=none] table[x=averageTimes,y=objLowerMedianVals] {2agent_ipmrandom.dat};
\addplot [fill=ipmrandom_color!10] fill between[of=upper and lower];


\addplot [name path=upper, draw=none] table[x=averageTimes,y=objUpperMedianVals] {2agent_proxhc.dat};
\addplot [name path=lower, draw=none] table[x=averageTimes,y=objLowerMedianVals] {2agent_proxhc.dat};
\addplot [fill=proxhc_color!10] fill between[of=upper and lower];

\draw[dashed, thick, color=black] (axis cs:2.2, 1e-1) -- (axis cs:2.2, 1e3);

\end{loglogaxis}
\end{tikzpicture}
\caption{Convergence of the objective function value (2 agents).}
\end{subfigure}

\bigskip

\begin{subfigure}{\linewidth}
\pgfplotsset{every tick label/.append style={font=\scriptsize}}
\begin{tikzpicture}
\begin{loglogaxis}[
    xlabel = \footnotesize Time (secs), ylabel=\footnotesize Constraints Violation,
    width=0.95\linewidth,
    height = 0.5\linewidth,
    xmin = 1e-1, xmax = 30
]
\definecolor{sqphc_color}{rgb}{0,0.4470,0.7410}
\addplot [mark=none, color=sqphc_color,line width=1pt] table[x=averageTimes,y=conMedianVals] {2agent_sqphc.dat};

\addplot [name path=upper,draw=none] table[x=averageTimes,y = conUpperMedianVals] {2agent_sqphc.dat};
\addplot [name path=lower,draw=none] table[x=averageTimes,y = conLowerMedianVals] {2agent_sqphc.dat};
\addplot [fill=sqphc_color!10] fill between[of=upper and lower];
\definecolor{sqprandom_color}{rgb}{0.4940,0.1840,0.5560}
\addplot [mark=none, color=sqprandom_color,line width=1pt] table[x=averageTimes,y=conMedianVals] {2agent_sqprandom.dat};

\addplot [name path=upper,draw=none] table[x=averageTimes,y = conUpperMedianVals] {2agent_sqprandom.dat};
\addplot [name path=lower,draw=none] table[x=averageTimes,y = conLowerMedianVals] {2agent_sqprandom.dat};
\addplot [fill=sqprandom_color!10] fill between[of=upper and lower];
\definecolor{ipmhc_color}{rgb}{0, 1, 0}
\addplot [mark=none, color=ipmhc_color,line width=1pt] table[x=averageTimes,y=conMedianVals] {2agent_ipmhc.dat};

\addplot [name path=upper,draw=none] table[x=averageTimes,y = conUpperMedianVals] {2agent_ipmhc.dat};
\addplot [name path=lower,draw=none] table[x=averageTimes,y = conLowerMedianVals] {2agent_ipmhc.dat};
\addplot [fill=ipmhc_color!10] fill between[of=upper and lower];

\definecolor{ipmrandom_color}{rgb}{0.9290,0.6940,0.1250}
\addplot [mark=none, color=ipmrandom_color,line width=1pt] table[x=averageTimes,y=conMedianVals] {2agent_ipmrandom.dat};

\addplot [name path=upper,draw=none] table[x=averageTimes,y = conUpperMedianVals] {2agent_ipmrandom.dat};
\addplot [name path=lower,draw=none] table[x=averageTimes,y = conLowerMedianVals] {2agent_ipmrandom.dat};
\addplot [fill=ipmrandom_color!10] fill between[of=upper and lower];


\definecolor{proxrandom_color}{rgb}{0,0,1}
\addplot [mark=none, color=proxrandom_color,line width=1pt] table[x=averageTimes,y=conMedianVals] {2agent_proxrandom.dat};

\addplot [name path=upper,draw=none] table[x=averageTimes,y = conUpperMedianVals] {2agent_proxrandom.dat};
\addplot [name path=lower,draw=none] table[x=averageTimes,y = conLowerMedianVals] {2agent_proxrandom.dat};
\addplot [fill=proxrandom_color!10] fill between[of=upper and lower];

\definecolor{proxhc_color}{rgb}{1,0,0}
\addplot [mark=none, color=proxhc_color,line width=1pt] table[x=averageTimes,y=conMedianVals] {2agent_proxhc.dat};

\addplot [name path=upper,draw=none] table[x=averageTimes,y = conUpperMedianVals] {2agent_proxhc.dat};
\addplot [name path=lower,draw=none] table[x=averageTimes,y = conLowerMedianVals] {2agent_proxhc.dat};
\addplot [fill=proxhc_color!10] fill between[of=upper and lower];

\draw[dashed, thick, color=black] (axis cs:2.2, 1e-20) -- (axis cs:2.2, 1e5);

\end{loglogaxis}
\end{tikzpicture}
\caption{Convergence of the constraints violation (2 agents).}
\end{subfigure}

\bigskip

\begin{subfigure}{\linewidth}
\pgfplotsset{every tick label/.append style={font=\scriptsize}}
\begin{tikzpicture}
\begin{loglogaxis}[
    xlabel = \footnotesize Time (secs),
    ylabel=\footnotesize Objective Function Value,
    width=0.95\linewidth,
    height = 0.5\linewidth,
    xmin = 1e-1, xmax=8e2, xshift = 0.3cm
]
\definecolor{sqphc_color}{rgb}{0,0.4470,0.7410}
\addplot [mark=none, color=sqphc_color,line width=1pt] table[x=averageTimes,y=objMedianVals] {6agent_sqphc.dat};

\addplot [name path=upper,draw=none] table[x=averageTimes,y = objUpperMedianVals] {6agent_sqphc.dat};
\addplot [name path=lower,draw=none] table[x=averageTimes,y = objLowerMedianVals] {6agent_sqphc.dat};
\addplot [fill=sqphc_color!10] fill between[of=upper and lower];

\definecolor{sqprandom_color}{rgb}{0.4940,0.1840,0.5560}
\addplot [mark=none, color=sqprandom_color,line width=1pt] table[x=averageTimes,y=objMedianVals] {6agent_sqprandom.dat};

\addplot [name path=upper,draw=none] table[x=averageTimes,y = objUpperMedianVals] {6agent_sqprandom.dat};
\addplot [name path=lower,draw=none] table[x=averageTimes,y = objLowerMedianVals] {6agent_sqprandom.dat};
\addplot [fill=sqprandom_color!10] fill between[of=upper and lower];

\definecolor{ipmhc_color}{rgb}{0, 1, 0}
\addplot [mark=none, color=ipmhc_color,line width=1pt] table[x=averageTimes,y=objMedianVals] {6agent_ipmhc.dat};

\addplot [name path=upper,draw=none] table[x=averageTimes,y = objUpperMedianVals] {6agent_ipmhc.dat};
\addplot [name path=lower,draw=none] table[x=averageTimes,y = objLowerMedianVals] {6agent_ipmhc.dat};
\addplot [fill=ipmhc_color!10] fill between[of=upper and lower];

\definecolor{ipmrandom_color}{rgb}{0.9290,0.6940,0.1250}
\addplot [mark=none, color=ipmrandom_color,line width=1pt] table[x=averageTimes,y=objMedianVals] {6agent_ipmrandom.dat};

\addplot [name path=upper,draw=none] table[x=averageTimes,y = objUpperMedianVals] {6agent_ipmrandom.dat};
\addplot [name path=lower,draw=none] table[x=averageTimes,y = objLowerMedianVals] {6agent_ipmrandom.dat};
\addplot [fill=ipmrandom_color!10] fill between[of=upper and lower];

\definecolor{proxrandom_color}{rgb}{0,0,1}
\addplot [mark=none, color=proxrandom_color,line width=1pt] table[x=averageTimes,y=objMedianVals] {6agent_proxrandom.dat};

\addplot [name path=upper,draw=none] table[x=averageTimes,y = objUpperMedianVals] {6agent_proxrandom.dat};
\addplot [name path=lower,draw=none] table[x=averageTimes,y = objLowerMedianVals] {6agent_proxrandom.dat};
\addplot [fill=proxrandom_color!10] fill between[of=upper and lower];

\definecolor{proxhc_color}{rgb}{1,0,0}
\addplot [mark=none, color=proxhc_color,line width=1pt] table[x=averageTimes,y=objMedianVals] {6agent_proxhc.dat};

\addplot [name path=upper,draw=none] table[x=averageTimes,y = objUpperMedianVals] {6agent_proxhc.dat};
\addplot [name path=lower,draw=none] table[x=averageTimes,y = objLowerMedianVals] {6agent_proxhc.dat};
\addplot [fill=proxhc_color!10] fill between[of=upper and lower];

\draw[dashed, thick, color=black] (axis cs:20, 1e-12) -- (axis cs:20, 1e4);

\end{loglogaxis}
\end{tikzpicture}
\caption{Convergence of the objective function value (6 agents).}
\end{subfigure}

\bigskip

\begin{subfigure}{\linewidth}
\pgfplotsset{every tick label/.append style={font=\scriptsize}}
\begin{tikzpicture}
\begin{loglogaxis}[
    xlabel = \footnotesize Time (secs),
    ylabel=\footnotesize Constraints Violation,
    width=0.95\linewidth,
    height = 0.5\linewidth,
    xmin = 1e-1, xmax=8e2
]
\definecolor{sqphc_color}{rgb}{0,0.4470,0.7410}
\addplot [mark=none, color=sqphc_color,line width=1pt] table[x=averageTimes,y=conMedianVals] {6agent_sqphc.dat};

\addplot [name path=upper,draw=none] table[x=averageTimes,y = conUpperMedianVals] {6agent_sqphc.dat};
\addplot [name path=lower,draw=none] table[x=averageTimes,y = conLowerMedianVals] {6agent_sqphc.dat};
\addplot [fill=sqphc_color!10] fill between[of=upper and lower];

\definecolor{sqprandom_color}{rgb}{0.4940,0.1840,0.5560}
\addplot [mark=none, color=sqprandom_color,line width=1pt] table[x=averageTimes,y=conMedianVals] {6agent_sqprandom.dat};

\addplot [name path=upper,draw=none] table[x=averageTimes,y = conUpperMedianVals] {6agent_sqprandom.dat};
\addplot [name path=lower,draw=none] table[x=averageTimes,y = conLowerMedianVals] {6agent_sqprandom.dat};
\addplot [fill=sqprandom_color!10] fill between[of=upper and lower];

\definecolor{ipmhc_color}{rgb}{0, 1, 0}
\addplot [mark=none, color=ipmhc_color,line width=1pt] table[x=averageTimes,y=conMedianVals] {6agent_ipmhc.dat};

\addplot [name path=upper,draw=none] table[x=averageTimes,y = conUpperMedianVals] {6agent_ipmhc.dat};
\addplot [name path=lower,draw=none] table[x=averageTimes,y = conLowerMedianVals] {6agent_ipmhc.dat};
\addplot [fill=ipmhc_color!10] fill between[of=upper and lower];

\definecolor{ipmrandom_color}{rgb}{0.9290,0.6940,0.1250}
\addplot [mark=none, color=ipmrandom_color,line width=1pt] table[x=averageTimes,y=conMedianVals] {6agent_ipmrandom.dat};

\addplot [name path=upper,draw=none] table[x=averageTimes,y = conUpperMedianVals] {6agent_ipmrandom.dat};
\addplot [name path=lower,draw=none] table[x=averageTimes,y = conLowerMedianVals] {6agent_ipmrandom.dat};
\addplot [fill=ipmrandom_color!10] fill between[of=upper and lower];

\definecolor{proxrandom_color}{rgb}{0,0,1}
\addplot [mark=none, color=proxrandom_color,line width=1pt] table[x=averageTimes,y=conMedianVals] {6agent_proxrandom.dat};

\addplot [name path=upper,draw=none] table[x=averageTimes,y = conUpperMedianVals] {6agent_proxrandom.dat};
\addplot [name path=lower,draw=none] table[x=averageTimes,y = conLowerMedianVals] {6agent_proxrandom.dat};
\addplot [fill=proxrandom_color!10] fill between[of=upper and lower];

\definecolor{proxhc_color}{rgb}{1,0,0}
\addplot [mark=none, color=proxhc_color,line width=1pt] table[x=averageTimes,y=conMedianVals] {6agent_proxhc.dat};

\addplot [name path=upper,draw=none] table[x=averageTimes,y = conUpperMedianVals] {6agent_proxhc.dat};
\addplot [name path=lower,draw=none] table[x=averageTimes,y = conLowerMedianVals] {6agent_proxhc.dat};
\addplot [fill=proxhc_color!10] fill between[of=upper and lower];

\draw[dashed, thick, color=black] (axis cs:20, 1e-13) -- (axis cs:20, 1e5);

\end{loglogaxis}
\end{tikzpicture}
\caption{Convergence of the constraints violation (6 agents).}
\end{subfigure}

\caption{Convergence of the prox-linear method and benchmark methods with filter-based warm-starting (FW) and random input initialization for 100 Monte Carlo simulations. In the figure, the black dashed vertical line represents the median time spent obtaining the warm-starting. The solid lines represent the median values of the simulations, while the shaded areas indicate the interquartile range, with the lower bound at the first quartile (0.25 quantile) and the upper bound at the third quartile (0.75 quantile).}
\label{fig: 100_convergence}
\end{figure}

We compare the convergence performance of our proposed method, the prox-linear method with filter-based warm-starting (FW), against benchmark methods such as SQP and IPM across 100 different random seeds in Fig.~\ref{fig: 100_convergence}. For the proposed method, we use the \texttt{OSQP}~\cite{osqp} solver to solve the linearized optimization problem, while for SQP and IPM, we use \texttt{fmincon} in \texttt{MATLAB}. The prox-linear method with random input initialization converges the fastest. However, it may get stuck in sub-optimal trajectories. The proposed method can reduce its objective value by up to 72\% when they both converge. Additionally, it reduces the objective function value for SQP and IPM with random input initialization by approximately 96\% within 3 seconds for the two-agent case, and by approximately 96\% and 98\% within 30 seconds for the six-agent case. Its convergence time is reduced by at least 91\% compared to SQP and IPM for the two-agent case, and by at least 96\% for the six-agent case. Furthermore, it achieves a better local optimum at terminal time, with up to a 70\% reduction in the objective function value compared to SQP, and up to a 77\% reduction compared to IPM. Warm-starting can also assist SQP and IPM in finding lower objective function values, though it may take longer to converge and achieve small constraint violations in large-scale problems. Therefore, the improvement is not as significant as it is with the prox-linear method.



\section{Conclusion}
We presented a first-order method with filter-based warm-starting to solve the nonconvex trajectory optimization problem. The approach combines constraint-aware particle filtering, hierarchical clustering, and the prox-linear method. We demonstrated the proposed method on a multi-agent system with linear dynamics and nonconvex inequality constraints. Compared with benchmark methods, the proposed method achieves both faster convergence and a lower objective function value.

This work still has some limitations. For example, we only considered linear dynamics, quadratic objective functions, and uniform discretization. In the future, we plan to extend the current work to more general problems beyond these frameworks, as well as further improve the computational time of the proposed method.




\bibliographystyle{IEEEtran}
\bibliography{IEEEabrv,reference}

\begin{thebibliography}{10}
\providecommand{\url}[1]{#1}
\csname url@rmstyle\endcsname
\providecommand{\newblock}{\relax}
\providecommand{\bibinfo}[2]{#2}
\providecommand\BIBentrySTDinterwordspacing{\spaceskip=0pt\relax}
\providecommand\BIBentryALTinterwordstretchfactor{4}
\providecommand\BIBentryALTinterwordspacing{\spaceskip=\fontdimen2\font plus
\BIBentryALTinterwordstretchfactor\fontdimen3\font minus \fontdimen4\font\relax}
\providecommand\BIBforeignlanguage[2]{{%
\expandafter\ifx\csname l@#1\endcsname\relax
\typeout{** WARNING: IEEEtran.bst: No hyphenation pattern has been}%
\typeout{** loaded for the language `#1'. Using the pattern for}%
\typeout{** the default language instead.}%
\else
\language=\csname l@#1\endcsname
\fi
#2}}

\bibitem{malyuta2021advances}
D.~Malyuta, Y.~Yu, P.~Elango, and B.~A{\c{c}}{\i}kme{\c{s}}e, ``Advances in trajectory optimization for space vehicle control,'' \emph{Annual Reviews in Control}, vol.~52, pp. 282--315, 2021.

\bibitem{malyuta2021convex}
D.~Malyuta, T.~P. Reynolds, M.~Szmuk, T.~Lew, R.~Bonalli, M.~Pavone, and B.~A{\c{c}}{\i}kme{\c{s}}e, ``Convex optimization for trajectory generation: A tutorial on generating dynamically feasible trajectories reliably and efficiently,'' \emph{IEEE Control Systems Magazine}, vol.~42, no.~5, pp. 40--113, 2022.

\bibitem{wang2024survey}
Z.~Wang, ``A survey on convex optimization for guidance and control of vehicular systems,'' \emph{Annual Reviews in Control}, vol.~57, p. 100957, 2024.

\bibitem{el1996formulation}
A.~El-Bakry, R.~A. Tapia, T.~Tsuchiya, and Y.~Zhang, ``On the formulation and theory of the {Newton} interior-point method for nonlinear programming,'' \emph{Journal of Optimization Theory and Applications}, vol.~89, no.~3, pp. 507--541, 1996.

\bibitem{byrd1999interior}
R.~H. Byrd, M.~E. Hribar, and J.~Nocedal, ``An interior point algorithm for large-scale nonlinear programming,'' \emph{SIAM Journal on Optimization}, vol.~9, no.~4, pp. 877--900, 1999.

\bibitem{lukvsan2004interior}
L.~Luk{\v{s}}an, C.~Matonoha, and J.~Vl{\v{c}}ek, ``Interior-point method for non-linear non-convex optimization,'' \emph{Numerical Linear Algebra with Applications}, vol.~11, no. 5-6, pp. 431--453, 2004.

\bibitem{boggs1995sequential}
P.~T. Boggs and J.~W. Tolle, ``Sequential quadratic programming,'' \emph{Acta Numerica}, vol.~4, pp. 1--51, 1995.

\bibitem{gill2011sequential}
P.~E. Gill and E.~Wong, ``Sequential quadratic programming methods,'' in \emph{Mixed Integer Nonlinear Programming}.\hskip 1em plus 0.5em minus 0.4em\relax New York, NY: Springer, 2011, pp. 147--224.

\bibitem{morgan2014model}
D.~Morgan, S.-J. Chung, and F.~Y. Hadaegh, ``Model predictive control of swarms of spacecraft using sequential convex programming,'' \emph{Journal of Guidance, Control, and Dynamics}, vol.~37, no.~6, pp. 1725--1740, 2014.

\bibitem{vitus2008tunnel}
M.~Vitus, V.~Pradeep, G.~Hoffmann, S.~Waslander, and C.~Tomlin, ``Tunnel-{MILP}: Path planning with sequential convex polytopes,'' in \emph{AIAA Guidance, Navigation and Control Conference and Exhibit}, 2008, p. 7132.

\bibitem{askari2021nonlinear}
I.~Askari, S.~Zeng, and H.~Fang, ``Nonlinear model predictive control based on constraint-aware particle filtering/smoothing,'' in \emph{American Control Conference}, 2021, pp. 3532--3537.

\bibitem{askari2022sampling}
I.~Askari, B.~Badnava, T.~Woodruff, S.~Zeng, and H.~Fang, ``Sampling-based nonlinear {MPC} of neural network dynamics with application to autonomous vehicle motion planning,'' in \emph{American Control Conference}, 2022, pp. 2084--2090.

\bibitem{askari2023model}
I.~Askari, X.~Tu, S.~Zeng, and H.~Fang, ``Model predictive inferential control of neural state-space models for autonomous vehicle motion planning,'' \emph{arXiv preprint arXiv:2310.08045 [cs.RO]}, 2023.

\bibitem{drusvyatskiy2018error}
D.~Drusvyatskiy and A.~S. Lewis, ``Error bounds, quadratic growth, and linear convergence of proximal methods,'' \emph{Mathematics of Operations Research}, vol.~43, no.~3, pp. 919--948, 2018.

\bibitem{drusvyatskiy2019efficiency}
D.~Drusvyatskiy and C.~Paquette, ``Efficiency of minimizing compositions of convex functions and smooth maps,'' \emph{Mathematical Programming}, vol. 178, pp. 503--558, 2019.

\bibitem{nielsen2016hierarchical}
F.~Nielsen, ``Hierarchical clustering,'' in \emph{Introduction to HPC with MPI for Data Science}, ser. Undergraduate Topics in Computer Science.\hskip 1em plus 0.5em minus 0.4em\relax Cham, Switzerland: Springer, 2016, pp. 195--211.

\bibitem{nocedal1999numerical}
J.~Nocedal and S.~J. Wright, \emph{Numerical Optimization}.\hskip 1em plus 0.5em minus 0.4em\relax New York, NY: Springer, 1999.

\bibitem{osqp}
B.~Stellato, G.~Banjac, P.~Goulart, A.~Bemporad, and S.~Boyd, ``{OSQP}: An operator splitting solver for quadratic programs,'' \emph{Mathematical Programming Computation}, vol.~12, no.~4, pp. 637--672, 2020.

\end{thebibliography}

\end{document}